\documentstyle{amsart}
\def\Q{{\bold Q}}
\def\Z{{\bold Z}}
\def\C{{\bold C}}
\def\R{{\bold R}}

\def\per{\mbox{per}}

\def\Aut{\mbox{Aut}}
\def\Hom{\mbox{Hom}}
\def\Lie{\mbox{Lie}}
\def\Supp{\mbox{Supp}}
\def\supp{\mbox{supp}}
\def\GL{\mbox{GL}}

\def\exp{\mbox{exp}}
\def\log{\mbox{log}}
\def\Log{\mbox{Log}}
\def\Alb{\mbox{Alb}}
\def\inv{\mbox{inv}}

\def\res{\mbox{res}}
\def\Lie{\mbox{Lie}}

\def\Div{\mbox{Div}}
\def\div{\mbox{div}}
\def \G{{\mathbf G}}
\def\Spec{\mbox{Spec}}

\title{p-adic abelian integrals and commutative Lie groups}

\author{Yu. G. Zarhin}
\address{The Pennsylvania State University, Department of Mathematics, 
University Park, PA 16802, USA
\newline
\indent Institute for Mathematical Problems in Biology, Russian Academy of Sciences,
Pushchino, Moscow Region, 142292, Russia}
\email{zarhin\char`\@math.psu.edu}
\dedicatory{To Yuri Ivanovich Manin}
\begin{document}

\maketitle

The aim of this paper is to propose an ``elementary" approach to Coleman's
theory of $p-$adic abelian integrals \cite{C1}, \cite{CG}. Our main tool is
a theory of commutative $p-$adic Lie groups (the logarithm map); we  use
neither {\sl dagger analysis} nor Monsky-Washnitzer cohomology theory. Notice
 that we also treat the case of bad reduction. 
We  discuss interrelations
between $p-$adic abelian integrals of of the third kind and N\'eron pairings on
abelian varieties.

A preliminary version of this paper appeared as a preprint \cite{Z4} in 1990.

{\sl Acknowledgments. } I am deeply grateful to Yu. I. Manin,
 D. Bertrand and A. N. Parshin for their interest in this paper.

\section{Logarithm maps}
Let $p$ be a prime, $\Q_p$ the field of $p-$adic numbers, $\C_p$ the completion
 of its algebraic closure. 
Let $K$ be a complete subfield of $\C_p$. Clearly, $K$ contains $\Q_p$. We will always deal with the valuation map 
$$v:K^*\to\Q$$ 
normalized by the condition $v(p)=1$. We will view $v$ as a homomorphism of the (multiplicative) $K-$Lie group $K^*$ into the discrete (additive) $K-$Lie group $\Q$ of rational numbers.

Let $G$ be a commutative $K-$algebraic group, $\Lie(G)$ its Lie algebra,
 $\Omega^1_{\inv}(G)$ the $K-$vector space of invariant differential forms of degree 1 on $G$. One may identify naturally $\Omega^1_{\inv}(G)$  with the dual of $\Lie(G)$, i.e.,
$$\Omega^1_{\inv}(G) =\Hom_K(\Lie(G),K).$$
Let us consider the $K-$analytic Lie group $G(K)$. Its Lie algebra coincides
 with $\Lie(G)$ (\cite{SH}, Ch. 2, Sect. 2.3 or \cite{Weil1}, Appendix III).
 If $u: G\to H$ is a homomorphism of $K-$algebraic groups, then it induces
 a homomorphism $G(K)\to H(K)$ of the $K-$analytic Lie groups and the corresponding algebraic and analytic tangent maps $du:\Lie(G)\to\Lie(H)$ coincide.

Recall (\cite{Bourbaki}, Chapitre III, 7.6) the properties of the logarithm map
$$\log_{G(K)}: G(K)_f\to \Lie(G).$$
Here $G(K)_f$ is the smallest open subgroup of $G(K)$ such that the quotient
 $G(K)/G(K)_f$ does not contain non-zero torsion elements,
 $\log_{G(K)}$ is a $K-$analytic homomorphism, whose tangent map
$$d\log_{G(K)}:\Lie(G)\to \Lie(\Lie(G))=\Lie(G)$$
is the identity map.  These properties determine $G(K)_f$ and $\log_{G(K)}$ uniquely. 

Here are some examples and remarks.

\begin{enumerate}
\item 
Let $G=\Spec K[t,t^{-1}]$ be the multiplicative group $\G_m$. Then
 $G(K)=K^*, \Lie(G)=K\cdot t \frac{d}{dt}\quad =K$ and $G(K)_f$ coincides with
 the group of units $U_K=\{x\in K^*\mid v(x)=0\}$ 
 (\cite{Bourbaki}, Ch. 3, Sect. 7, Ex. 3). The logarithm map
 $\log_{K^*}: U_K\to K$ coincides with the usual logarithm $\log$
 (defined via the usual convergent power series) on the subgroup of 
 ``principal" units $U_K^1=\{x\in K^*\mid v(1-x)>0\}$.
\item
Let $G=\Spec\ K[t]$ be the additive group $\G_a$ then $G(K)=G_f(K)=K$, $\Lie(G)=K\cdot\frac{d}{dt}\quad =K$ and $\log_K  : K\to K$ is the identity map.
\item
Let $G$ be an abelian variety. Then $G(K)=G(K)_f$, because if $U$ is an open
 subgroup of $G(K)$ then the quotient $G(K)/U$ is a torsion group 
(see Coleman(\cite{C2}; the proof is based on results of Raynaud \cite{R}). If
 $K$ is a finite algebraic extension of $\Q_p$, i.e., is locally compact, the equality follows easily from the compactness of $G(K)$.
\item
Let $G,H$ be commutative $K-$algebraic groups. Then
 $(G\times H)(K)=G(K)\times H(K)$, $\Lie(G\times H)=\Lie(G)\times\Lie(H),
  (G\times H)(K)_f=G(K)_f\times H(K)_f$ and in obvious notation
 $\log_{(G\times H)(K)}=(\log_{G(K)},\log_{H(K)})$.
\item
If $u: G\to H$ is a homomorphism of $K-$algebraic groups, then $u(G(K)_f)\subset H(K)_f$ and 
$$du\ \log_{G(K)}=\log_{H(K)} u.$$
\item
If $L$ is a closed algebraic $K-$subgroup of $G$ then 
$$L(K)_f=G(K)_f\bigcap L(K)$$ 
and the restriction of $\log_{G(K)}$ to $L(K)_f$ coincides with
$$\log_{L(K)}:L(K)_f\to \Lie(L)\subset \Lie(G).$$
\item
Let  $K' \subset \C_p$ be a complete overfield of $K$ and $G(K')$ be the $K'-$analytic Lie group of all $K'-$points of $G$. Then
$$G(K)_f=G(K')_f\bigcap G(K)$$ 
and the restriction of $\log_{G(K')}$ to $G(K)_f$ coincides with
$$\log_{G(K)}:G(K)_f\to \Lie(G)\subset \Lie(G)\otimes_K K'=\Lie(G(K')).$$

\end{enumerate}

Since $G(K)/G(K)_f$ is torsion-free, one may easily extend $\log_{G(K)}$ to a homomomorphism
$G(K)\to\Lie(G)$, which automatically is analytic and whose tangent map is the identity map. Now we explain how this extension can be chosen canonically for all $G$ , in order to keep the functoriality properties.

\vskip 1cm

{\bf Step 1}  We start with the case of $G=\G_m$. Notice that $v$ defines an isomorphism between $C^*/U_K$ and the additive group $\Q$. So, in order to extend $\log$ to a homomorphism $\C_p^*\to\C_p$, which sends $K^*$ into $K$, one has only to
choose $c \in K$ and  fix a branch of $p-$adic logarithm
$$\log=\log^{(c)}:\C_p^*\to \C_p,$$
by putting
$$  \log(p)=\log^{(c)}(p):=c$$
\cite{C2}, \cite{Z1}.   
(If $\Log:K^*\to K$ is another branch of the logarithm then
$$\Log=\log+(\Log(p)-\log(p))v .)$$ 
If $K' \subset \C_p$ is a complete overfield of $K$ then $\log({K'}^*)\subset K'$. Notice, that for any automorphism $\sigma$ of $K'/K$  which preserves the absolute value (i.e., $v$) ,
$$\log(\sigma(x))=\sigma((\log(x)) \quad\mbox{\rm for all }  x \in {K'}^*.$$

\vskip 1cm

 {\bf Step 2} Let us extend the logarithm map for a split torus $G=\G_m^r$. We
 have $G(K)=(K^*)^r, \Lie(G(K))=K^r$ and extend the logarithm map as follows.
$$\log^{(c)}_{(K^*)^r}(x_1, \ldots , x_r)=(\log(x_1), \ldots \log(x_r) ).$$
One may easily check that if $u:G=\G_m^r\to H=\G_m^n$ is a homomorphism of
 algebraic $K-$tori then 
$$\log^{(c)}_{ H(K)_f}(ux)=du(\log^{(c)}_{G(K)}(x)).$$
Indeed, there exists an integral $r\times n$-matrix $(a_{ij})$ such that $u$ acts via
$$((x_1,\ldots ,x_r)\mapsto (y_1,\ldots ,y_n), \quad  y_j=\prod_{i=1}^r x_i^{a_{ij} }\quad (j=1,\ldots n)$$
and the matrix of the tangent linear map 
$$du: K^r=\Lie(G)\to\Lie(H)=K^r$$
coincides with $(a_{ij})$. Now, one has only to notice that
$$\log^{(c)}(\prod_{i=1}^r x_i^{a_{ij}})=\sum_{i=1}^r a_{ij} \log^{(c)}(x_i ).$$

Let  $K' \subset \C_p$ be a complete overfield of $K$. Then we may consider the
 $K'-$analytic Lie group $G(K')=({K'}^*)^r$ with the $K'-$Lie algebra $\Lie(G)_{K'}
={K'}^r$ and define  the extended logarithm map
$$\log^{(c)}_{({K'}^*)^r}: G(K')=({K'}^*)^r\to\Lie(G)_{K'}={K'}^r, $$
$$ (x_1, \ldots , x_r)\mapsto(\log(x_1), \ldots \log(x_r) ).$$
Clearly, $\log^{(c)}_{({K'}^*)^r}$ coincides with $\log^{(c)}_{(K^*)^r}$ on
 $G(K)=(K^*)^r\subset ({K'}^*)^r=G(K')$.

Notice that for any automorphism $\sigma$ of $K'/K$  which preserves the absolute value (i.e., $v$),
$$\log^{(c)}_{({K'}^*)^r}(\sigma(x))=\sigma(\log^{(c)}_{({K'}^*)^r}(x))
 \quad\mbox{\rm for all }  x \in ({K'}^*)^r=G(K').$$

\vskip 1cm

{\bf Step 3} Now, we do the case of arbitrary algebraic $K-$torus $T$. Let us
 choose a fiinite Galois extension $K'/K$ , which sits in $\C_p$ and such that
 $T$ splits over $K'$, i.e., there exists an isomorphism of $K'-$algebraic tori
 $\phi:T_{K'}\cong \G_m^r$. This allows us to identify $T(K')$ with $({K'}^*)^r$ and
 $\Lie(T)_{K'}=\Lie(T)\otimes_K K'$ with ${K'}^r$. Then there exists a
 homomorphism (cocycle) of the Galois group $\Gamma$ of $K'/K$ 
$$\Gamma\to \Aut(\G_m^r)=\GL(r,\Z),\quad \sigma\mapsto \psi_{\sigma},$$
 such that the natural Galois action on $T(K')$ is defined by the formula
$$\psi_T(\sigma):x\mapsto \psi_{\sigma}(\sigma(x)) \quad
 \mbox{for all } \sigma \in \Gamma, x \in ({K'}^*)^r = T(K').$$
(Here $\sigma(x)=\sigma(x_1,\ldots , x_r)=(\sigma x_1, \ldots , \sigma x_r)$
 for  $x=(x_1,\ldots , x_r) \in ({K'}^*)^r$.)
Clearly, 
$$T(K)=\{x\in ({K'}^*)^r\mid x=\psi_{\sigma}(\sigma(x))
 \forall \sigma \in \Gamma\},$$
$$\Lie(T)=\{b \in {K'}^r\mid  b=\psi_{\sigma}(\sigma(b))
\forall \sigma \in \Gamma\}.$$
It follows easily that
$$\log^{(c)}_{({K'}^*)^r}(T(K)) \subset \Lie(T)$$
and therefore one may define
 $$\log^{(c)}_{T(K)}:T(K)\to\Lie(T)$$
 as the restriction of $\log^{(c)}_{({K'}^*)^r}$ to $T(K)$. One may easily
 check, using the results of the previous step, that the definition of
 $\log^{(c)}_{T(K)}$ does not depend on the choice of isomorphism
 between $T_{K'}$ and $\G_m^r$.

\vskip 1cm
{\bf Step 4} Now, assume that $G$ is a connected commutative linear algebraic group. Then $G$ is isomorphic to the product $T\times \G_a^n$ where $T$ is algebraic $K-$torus. Then
$$G(K)=T(K)\times K^n, \quad \Lie(G)=\Lie(T)\times K^n$$
and we define
$$\log^{(c)}_{T(K)\times K^n}:T(K)\times K^n\to \Lie(T)\times K^n$$
by the formula
$$\log^{(c)}_{T(K)\times K^n}(x, a_1, \ldots , a_n)=(\log^{(c)}_{T(K)}(x), a_1, \ldots , a_n).$$

\vskip 1cm

{\bf Step 5} Let $G$ be a connected commutative algebraic $K-$group. Then, by a theorem of Chevalley, $G$ sits in a short exact sequence
$$0\to L\to G\to A\to 0$$
where $L$ is a connected linear algebraic group and $A$ is an abelian variety. This gives rise to the short exact sequence of the $K-$Lie algebras
$$0\to\Lie(L)\to\Lie(G)\to \Lie(A)\to 0$$
and to the exact sequence of the $K-$Lie groups
$$0\to L(K)\to G(K)\to A(K).$$
It follows easily that the image $U$ of $G(K)_f$ in $A(K)$ is an open subgroup of $A(K)$. Recall that the quotient  $A(K)/U$ is a torsion group. This implies that the quotient $G(K)/(L(K)\ G(K)_f)$ is also a torsion group, because it is isomorphic to a subgroup of  $A(K)/U$.
It follows that for each $x\in G(K)$ there exist a positive integer $m$ and an element $x_f\in G(K)_f$ such that 
$$y=x^m (x_f)^{-1} \in L(K).$$
In other words,
$$x^m =x_f y, \quad x_f \in G_f(K), y\in L(K).$$
Now, let us put
$$\log^{(c)}_{G(K)}(x)=\frac{1}{m}(\log_{G(K)}(x_f) +\log^{(c)}_{L(K)}(y) )\in \Lie(G).$$
If
$$x^n =x'_f y', \quad x_f \in G_f(K), y\in L(K)$$
then $x_f^n y^n=(x'_f)^m (y')^m$, i.e.,
$$z=x_f^n (x'_f)^{-m} =(y')^m y^{-n} \in G(K)_f \bigcap L(K)=L(K)_f$$
and therefore
$$n\cdot \log_{G(K)}(x_f)-m\cdot \log_{G(K)}(x'_f)=\log_{G(K)}(z)=\log_{L(K)}(z)=\log^{(c)}_{L(K)}(z).$$
This implies that
$$\frac{1}{m}(\log_{G(K)}(x_f) +\log^{(c)}_{L(K)}(y) )=\frac{1}{n}(\log_{G(K)}(x'_f) +\log^{(c)}_{L(K)}(y') )\in \Lie(G),$$
i.e., the definition of $\log^{(c)}_{G(K)}(x)$ does not depend on the choice
 of  $m$ and $x_f$. It also follows easily that $\log^{(c)}_{G(K)}$ is a
 homomorphism from $G(K)$ to $\Lie(G)$, which coincides with $\log_{G(K)}$ on $G(K)_f$.

\vskip 1cm

{\bf Step 6} Assume that $G$ is not necessarily connected commutative algebraic group. Let $G^0$ be its connected identity component. We write $n$ for the index of $G^0$ in $G$. Then
$\Lie(G)=\Lie(G^0)$ and $x^n\in G^0(K) \forall x\in G(K)$. We put

$$\log^{(c)}_{G(K)}(x)=\frac{1}{n} \log^{(c)}_{G^0(K)}(x^n).$$

\vskip 1cm

The following statement summarizes the results of the steps 1-6.

\vskip 1 cm

{\bf Theorem.} {\sl To each commutative} $K-${\sl algebraic group} $G$ {\sl one may attach a} $K-${\sl analytic homomomorphism}
$$\log^{(c)}_{G(K)}\to \Lie(G), $$
{\sl enjoying the following properties}:

\begin{enumerate}
\item
{\bf (Choice of branch)} {\sl If} $K=\G_m$ {\sl then} $G(K)=K^*$ {\sl and}
$$\log^{(c)}_{G(K)}=\log^{(c)}:K^*\to K=K\cdot t \frac{d}{dt}.$$
\item
{\bf (Logarithm property)} {\sl The restriction of} $\log^{(c)}_{G(K)}$ {\sl to }$G(K)_f$ {\sl coincides with} $\log_{G(K)}$. {\sl In particular, its tangent map coincides with the identity map} $\Lie(G)\to\Lie(G)$.
\item
{\bf (Functoriality)} {\sl If} $G\to H$ {\sl is a homomorphism of commutative} $K-${\sl algebraic groups then}
$$du \ \log^{(c)}_{G(K)}=\log^{(c)}_{H(K)} u. $$
\item
{\bf (Product formula)} {\sl If} $G, H$ {\sl are commutative} $K-${\sl algebraic groups then} $(G\times H)(K)=G(K)\times H(K), \Lie(G\times H)=\Lie(G)\times\Lie(H)$ {\sl and in obvious notation}
$$\log^{(c)}_{(G\times H)(K)}=(\log^{(c)}_{G(K)},\log^{(c)}_{H(K)}).$$
\item
{\bf (Ground field extension)} {\sl Let}  $K' \subset \C_p$ {\sl be a complete overfield of} $K$ {\sl and} $G'=G\times K'$ {\sl be the corresponding commutative} $K'-${\sl algebraic group. Then} $G'(K')=G(K)$ is a commutative $K'-${\sl analytic Lie group}, {\sl the} $K'-${\sl Lie algebra} $\Lie(G')=\Lie(G)\otimes_K  K'$ {\sl and the restriction of} $\log^{(c)}_{G'(K')}$ to $G(K)\subset G(K')=G'(K')$ {\sl coincides with}
$$\log^{(c)}_{G(K)}:G(K)\to \Lie(G)\subset \Lie(G)\otimes_K  K'=\Lie(G'). $$
{\sl In addition, if} $\sigma$ {\sl is an automorphism of the field extension} $K'/K$, {\sl which preserves the absolute value,  then}
$$\log^{(c)}_{G'(K')}(\sigma(x))=\sigma(\log^{(c)}_{G'(K')}(x)) \quad\mbox{\rm for all }  x \in G(K')=G'(K').$$
\end{enumerate}

{\sl These properties determine the homomorphisms}
$\log^{(c)}_{G(K)}$ {\sl uniquely}.

\section{N\'eron pairings as periods of logarithm maps}

Let us consider an extension $G$ of an abelian variety $A$ by the multiplicative group $\G_m$, i.e., assume that $G$ sits in a short exact sequence
$$0\to\G_m\to G\to A\to 0.$$
The Hilbert's theorem 90 implies the exactness of the corresponding sequence
$$0\to K^*\to G(K)\to A(K)\to 0.$$
Recall that the corresponding exact sequence of the $K-$Lie algebras
$$0\to K\ \to \Lie(G) \to \Lie(A)\to 0$$
is also exact. Let us consider the difference of the two branches of the logarithm map
$$\per_G:\log^{(c+1)}_{G(K)}-\log^{(c)}_{G(K)}: G(K)\to \Lie(G)$$
attached to $c+1$ and $c$ respectively. Clearly, $\per_G$ is a locally constant homomorphism, which kills $G(K)_f$ and coincides with 
$$v:K^*\to \Q\subset K\subset\Lie(G)$$
 on $K^*\subset G(K)$. Since $G(K)/(K^*\ G(K)_f)$ is a torsion group (see step 5 of the previous section),
$$\per_G(G(K))\subset \Q\subset K\subset \Lie(G)$$
and does not depend on the choice of $c$. Indeed, if $x\in G(K)$ and a positive integer $m$, $a\in K^*$ and $x_f\in G(K)_f$ satisfy $x^m=a x_f \in G(K)$ then
$$m\per_G(x)=\per_G(x_f)+\per(a)=v(a)\in \Q$$
and therefore
$$\per_G(x)= v(a)/m \in \Q$$
does not depend on the choice of $c$.

Since $\per_G$ is locally constant, it may be viewed as a continuous homomorphism
$$\per_G: G(K)\to \Q \subset \R,$$
whose restriction to $K^*\subset G(K)$ coincides with
$$v:K^*\to\Q\subset \R.$$
These properties determine $\per_G$ uniquely. Indeed, if $\per':G(K)\to\R$ is
 a a locally constant homomorphism, coinciding with $v$ on $K^*$ then the
 difference $\per_G-\per':G(K)\to\R$ kills $K^*$ and a certain open subgroup
 $U'$ of $G(K)$. The image $U$ of $U'$ in $A(K)$ is an open subgroup and
 therefore $A(K)/U$ is a torsion group. The quotient $G(K)/(K^*\ U')$ is
 isomorphic to a subgroup of $A(K)/U$ and therefore is also a torsion group. Since $\per_G-\per'$ kills $K^*\ U'$, one has only to notice that $\R$ is uniquely divisible and therefore $(\per_G-\per')(G(K))=0.$

Recall \cite{S3}(see also \cite{L}) that there exists a linear bundle $L$ over
 $A$, which is algebraically equivalent to zero and such that $G$, viewed as a
principal $\G_m-$bundle over $A$, 
 is  isomorphic to the principal $\G_m-$bundle
 $$L^*=L\setminus \{\mbox{\rm the image of zero section}\}\to A.$$
  Further, we will identify $G$ with $L^*$.

Let us consider the absolute value on $K$ defined by the formula
$$\mid a\mid _p = \exp (-v(a)) \forall x \in K^*.$$

 There is a {\sl canonical local height} (\cite{Z3}, Sect. 3; see also \cite{Z4})
$$\hat h_{-v,L}:L^*(K)=G(K)\to \R,$$
which is a continuous homomorphism, coinciding with $-v$ on $K^*$.  Since each neighborhood of the identity element of $G(K)$ contains an open subgroup,
 $\hat h_{-v,L}$ must kill a certain open subgroup of $G(K)$, i.e., $\hat h_{-v,L}$ is locally constant. This implies that
              $$\hat h_{-v,L}=-\per_G.$$
Warning: $v$ as  defined in \cite{Z3} coincides with our $-v$!

 Recall \cite{Z1} (see also \cite{L}, \cite{Z3} ) a construction of N\'eron
 pairings $\langle D,\ {\frak a}\rangle_{-v}$ \cite{N} in terms of $\hat h_{-v,L}$.
 Here $D$  is a divisor on $X$, algebraically equivalent to zero, whose rational equivalence class coincides with (isomorphism class of)
 $L$, and ${\frak a} = \sum a_xx$  is a $0-$dimensional cycle of degree zero on $X$, whose support lies in $X(K)$ and does not meet $\Supp(D)$.  Let $s_D$ be a non-zero rational section of
$L$, whose divisor coincides with $D$. Recall that $s_D$  is determined uniquely up to multiplication by an element of $K^*$. Then
$$\langle D,\ {\frak a}\rangle_{-v}= -\hat h_{-v,L}\bigg(\sum a_x s_D(x)\bigg)=-\sum 
a_x \hat h_{-v,L}(s_D(x)).$$ 
It follows that
$$\langle D,\ {\frak a}\rangle_{-v}=\sum 
a_x \per_G(s_D(x)) \in \Q.$$
The rationality of  $\langle D,\ {\frak a}\rangle_{-v}$ is well-known \cite{N} (see also \cite{L}) when $K$ is a discrete valuation field.

\section{Differentials of the third kind}
Let $V$ be an absolutely irreducible projective variety over $K$. We assume
 that the set $V(K)$ of $K-$rational points is non-empty. Clearly, $V(K)$ is
 dense in $V$ in the Zariski topology. Let $\Div(V)$ be the group of divisors on $V$.
Recall that each divisor $D\subset \Div(V)$ can be uniquely presented as a formal
linear combination $\sum c_Z Z$ where $Z$ runs through the set $V^{(1)}$ of
 $K-$irreducible
closed subvarieties of codimension 1, $c_Z$ are integers vanishing for all but
finitely many $Z$. We will mostly deal with the tensor product
 $\Div_A(V)\otimes K$, whose elements are formal  linear combinations
 $\sum c_Z Z$ where $c_Z$ are elements of $K$, vanishing for all but finitely
many $Z\in V^{(1)}$. To each non-zero
 $D=\sum_{Z\in V^{(1)}}c_Z Z\in\Div (V)\otimes K$
we attach the finite set
$\mbox{supp} (D)=\{Z\in V^{(1)}\mid c_Z\ne 0\}\subset V^{(1)}$
and the codimension 1 subvariety $\mbox{Supp}(D)=\bigcup_{Z\in\supp(D)} D \subset V$.
Recall \cite{S2}, \cite{FW} that a K\"ahler
 differential $\omega$ on $V$ is a {\sl differential of the third kind} if locally
 in the Zariski topology one may represent $\omega$ as a sum
$$\omega=\omega_{\mbox{reg}} + \sum c_i df_i/f_i $$
of regular differential $\omega_{\mbox{reg}}$ and a finite linear combination of
logarithmic differentials $df_i/f_i$ of non-zero rational functions 
$f_i\in K(V)$ with coefficients $c_i\in K$. We write $\res(\omega)$ for the
residue of $\omega$ \cite{S1,S2,FW}; by definition
$$\res(\omega) \in \Div_a(V) \otimes K \subset \Div(V)\otimes K$$
where $\Div_a(V)$ is the group of divisors on $V$, algebraically equivalent to
zero. The differential $\omega$ is regular outside $\Supp(\res(\omega))$. For
example, if $\omega=df/f$ is the logarithmic differential of (non-constant)
rational function $f$, then $\res(\omega)=\div(f)$ is the divisor of $f$ and
$\omega$ is regular outside $\Supp(\div(f))$. As usual, we agree that empty set is the support
 of zero divisor and zero differential.

 Further, we assume that the set $V(K)$ is non-empty. Let $S\subset V^{(1)}$ be a finite set
 of irreducible codimension 1 subvarieties in $V$. We write $K[S]$ for the finite-dimensional
  $K-$vector subspace of $\Div_A(V)\otimes K$, generated by elements of $S$, i.e., the subspace,
   consisting of elements, whose $\supp$ is contained in $S$. We write $\Omega_{3,S}(V)$ for the
   finite-dimensional $K-$vector space of differentials $\omega$ of the third kind on $V$ with
   $\res(\omega)\in K[S]$. We write $K[S]_a$ for the intersection of $k[S]$ and
   $\Div_a(V)\otimes K$ in $\Div(V)\otimes K$. Clearly, $\Omega_{3,S}(V)$ contains the space
   $\Omega^1(V)$ of differentials of the first kind, i.e., everywhere regular differentials on
   $V$ and $\res(\Omega_{3,S}(V))\subset K[S]_a$.

It follows from \cite{S1}, Lemma 6 and its proof (see also \cite{S2} and 
\cite{FW}, especially, pp. 197--201) that there exists an extension
$$ 0 \to T \to G_S \to \Alb (V) \to 0 $$
$G_S$ of the Albanese variety $\Alb (V)$ of $V$ by a $K-$torus $T$ and a rational map $f: V \to G_S$, which is regular on the complement of union of $D \in S$ and satisfies the following properties. If $\Omega^1_{\inv}(G_S)$ is the $K-$vector space of invariant differential forms of degree 1 on $G_S$ then $f^*\Omega^1_{\inv}(G_S)=\Omega_{3,S}(V)$ and the map $\omega\mapsto f^*\omega$ is an isomorphism between $\Omega^1_{\inv}(G_S)$ and $\Omega_{3,S}(V)$. In other words, for each $\omega\in \Omega_{3,S}(V)$ there exists exactly one invariant form $\omega_{\inv}\in\Omega^1_{\inv}(G_S)$ with $f^*\omega_{\inv}=\omega$. The rational map $f$ is unique, up to a translation by an element of $G_S(K)$. The   algebraic group $G_S$ is unique, up to an isomorphism.

So, in order to integrate $\omega\in \Omega_{3,S}(V)$ it suffices to know how to integrate $\omega_{\inv}$ on $G_S$ and just put
$$\int_P^Q \omega =\int_{f(P)}^{f(Q)}\omega_{\inv} \in K.$$

Now, let us use the duality between $\Lie(G_S)$ and $\Omega^1_{\inv}(G_S)$. Let us choose a branch $\log^{(c)}$ of the $p-$adic logarithm and put
$$\int_x^y \omega_{\inv} := \omega_{\inv}(\log^{(c)}_{G_S(K)}(y)-\log^{(c)}_{G_S(K)}(x)) \in K.$$

The $K-$valued function(s) $u(P,Q)=u(P, Q,\omega) =\int_{f(P)}^{f(Q)}\omega_{\inv}$ obviously enjoys the following natural properties:

\begin{itemize}
\item (Newton--Leibnitz rule) $u(P,Q)$ is a $K-$analytic function in $Q$ and its differential coincides with $\omega$;

\item ($K-$linearity) $u(P,Q)=u(P,Q,\omega)$ is a $K-$linear function in $\omega$;
\item (Additivity) $u(P,Q)+u(Q,R)=u(P,R)$ for all $P,Q,R \in V(K)\setminus \Supp(\res(\omega))$.
\item (Change of variables) Let $u: V\to W$ be a regular dominant map of absolutely irreducible projective
algebraic $K-$varieties, $\omega'$ a differential of the third kind on $W$ such that $u^*\omega'=\omega$.
Then
$$\int_P^Q\omega=\int_{u(P)}^{u(Q)}\omega' $$
whenever  both integrals are defined.
\item
(Choice of branch) If $\omega=dz/z$ is a logarithmic differential of a rational function $z$ on $V$ then
$$\int_P^Q dz/z= \log^{(c)}(z(Q))-\log^{(c)}(z(P)).$$
\end{itemize}

\vskip 1cm {\bf Remark.} Using Lie theory over the fields of real and complex numbers (maximal compact subgroups), Bertrand and Philippon  have constructed differentials of the third kind with prescribed residues and purely imaginary periods (\cite{BP}, Sect. 4, especially Remarque 5).

\section{Closed forms}
Let $X$ be an absolutely irreducible smooth (not necessarily projective)
 variety over $K$. We assume that the set $X(K)$ is non-empty. Let $\omega$ be
a {\sl closed} regular differential form of degree 1 on $X$. In this section we define
 $\int_P^Q \omega$ for $P,Q\in X(K)$. According to  Satz 5 of
 \cite{FW}, p. 197, there exist a commutative algebraic $K-$group $G$,
 an  invariant form $\omega_{\inv}\in \Omega_{\inv}^1(G)$ and a regular map $f: X\to G$ such that
$$f^* \omega_{\inv}=\omega.$$
Now, we put
$$\int_P^Q \omega=\omega_{\inv}(\log^{(c)}_{G(K)}(f(Q))-\log^{(c)}_{G(K)}(f(P))) \in K.$$
One may easily deduce from  Satz 6 of \cite{FW}, p. 197 that
 $\int_P^Q \omega$ does not depend on the choice of $G, f$ and $\omega_{\inv}$. It
 also follows easily that the integral enjoys the Newton--Leibnitz, $K-$linearity, additivity, choice of branch and change of variables properties.


\begin{thebibliography}{99}
\bibitem{BP}
Bertrand D. and Philippon P. {\sl Sous-groupes alg\'ebriques  de groupes algebriques commutatifs}. Illinois J. Math. {\bf 32} (1988), 263--280.
\bibitem{Bourbaki}
Bourbaki N. Groupes et Alg\'ebres de Lie, Hermann, Paris, 1972.
\bibitem{C1}
Coleman R. F. {\em Torsion points on curves and $p-$adic abelian integrals}.
Ann. of Math. {\bf 121} (1985), 111--168.
\bibitem{C2}
Coleman R. F. {\em Reciprocity laws on curves}. Compositio Math. {\bf 72} (1989), 205--235.
\bibitem{CG}
Coleman R. F. and Gross B. $p$-{\sl adic heights on curves}. Advanced Studies
 in Pure Mathematics {\bf 17} (1989), 73--81.

\bibitem{FW}
Faltings G. and W\"ustholz G. {\em Einbettungen commutativer algebraischer
Gruppen und einige ihrer Eigenschaften}. J. reine angew. Math. {\bf 354} (1984),
175--205.

\bibitem{L}
Lang S. Fundamentals of Diophantine Geometry. Springer-Verlag, 1983.
\bibitem{N}
N\'eron A. {\em Quasi-functions et hauteurs sur les vari\'et\'es ab\'eliennes}.
Ann.of Math. {\bf 82} (1965), 249--331.
\bibitem{R}
Raynaud M. {\em Vari\'et\'es ab\'eliennes et g\'eom\'etrie rigide}. Actes Congres Intern.
Math. Tome 1 (1970), 473--477.
\bibitem{S1}
Serre J.-P. {\em Morphismes universels et vari\'et\'e d'Albanese}. S\'eminaire
C. Chevalley (1958/59), E.N.S.,Paris.
\bibitem{S2}
Serre J. P. {\em Morphismes universels et diff\'erentielles de troisi\`eme
 esp\'ece}.
S\'eminaire C. Chevalley (1958/59), E.N.S., Paris.
\bibitem{S3}
Serre J.P. Groupes alg\'ebriques et corps de classes. Paris, Hermann, 1959.

\bibitem{SH}
Shafarevich I.\ R. Basic Algebraic Geometry,
Springer-Verlag, Berlin-Heidelberg-New York, 1974 (= 
Osnovy algebraicheskoi geometrii, Nauka, Moscow, 1972).

\bibitem{Weil1}
 Weil A.\ Foundations of Algebraic Geometry,
AMS Colloquium Publications {\bf 29}, Revised Edition,
American Mathematical Society, Providence, R.\ I., 1962.

\bibitem{Z1}
Zarhin Yu. G. $p-${\em adic heights on abelian varieties}. S\'eminaire de
Th\'eorie des Nombres, Paris 1987-88 (C. Goldstein, ed.). Progress in Math.
 {\bf 81} (1990), 317--341.
 Birkh\"auser, Basel.

\bibitem{Z2}
Zarhin Yu. G. {\em N\'eron pairing and quasicharacters}. Izv. Akad. Nauk SSSR
ser. matem. {\bf 36} (1972), 497--509; English translation in Math. USSR Izv.
{\bf 6} (1972), 491--503.
\bibitem{Z3}
Zarhin Yu. G. {\em Local heights and N\'eron pairings}. Trudy Steklov Math.
 Inst. {\bf 208} (1995), 111-127.
\bibitem{Z4}
Zarhin Yu. G. {\em Local heights and abelian integrals}. Expos\'e  9 dans
``Problemes Diophantiens 88-89" (D. Bertrand, M. Waldschmidt), Publ. Math.
Univ. Paris VI {\bf 90} (1990), 15 pp.
\end{thebibliography}
\end{document}